\documentclass[12pt]{article}
\usepackage{amsmath,amssymb,amsthm,morefloats}
\usepackage{epsfig,psfrag,url}
\usepackage{graphicx}
\usepackage{verbatim}
\usepackage{subfigure}
\usepackage{somacros}

\numberwithin{equation}{section}

\def\Mobius{M\"obius\ }

\newcommand{\RHP}{RHP}
\def\CC{\mathcal C}
\newcommand{\goto}{\rightarrow}
\newcommand{\bigo}{{\mathcal O}}

\def\XXint#1#2#3{{\setbox0=\hbox{$#1{#2#3}{\int}$}
     \vcenter{\hbox{$#2#3$}}\kern-.5\wd0}}

\DeclareMathOperator{\diag}{diag}

\newenvironment{choices}{\left\{ \begin{array}{ll}}{\end{array}\right.}
\newcommand\when{&\text{if~}}
\newcommand\otherwise{&\text{otherwise}}

\newenvironment{mat}{\left(\begin{array}{ccccccccccccccc}}{\end{array}\right)}
\newcommand\bcm{\begin{mat}}
\newcommand\ecm{\end{mat}}

\newcommand{\bea}{\begin{eqnarray}}
\newcommand{\eea}{\end{eqnarray}}
\newcommand{\bean}{\begin{eqnarray*}}
\newcommand{\eean}{\end{eqnarray*}}
\newcommand{\ba}{\begin{array}}
\newcommand{\ea}{\end{array}}
\newcommand{\beqs}{\begin{equation*}\begin{split}}

\newtheorem{problem}{Problem}[section]

\newtheorem{algorithm}{Algorithm}[section]

\theoremstyle{definition}

\newtheorem{remark}{Remark}[section]
\newtheorem{definition}{Definition}[section]

\long\def\symbolfootnote[#1]#2{\begingroup%
\def\thefootnote{\fnsymbol{footnote}}\footnote[#1]{#2}\endgroup}

\begin{document}
\title{A Riemann--Hilbert approach to Jacobi operators and Gaussian quadrature}
\author{Thomas Trogdon\\
{\tt trogdon@cims.nyu.edu}\\
Courant Institute of Mathematical Sciences\\
 New York University\\
251 Mercer St.\\
New York, NY, 10012, USA \\
\phantom{.}\\
Sheehan Olver\\
{\tt Sheehan.Olver@sydney.edu.au}\\
School of Mathematics and Statistics\\
 The University of Sydney\\
NSW 2006, Australia}
\maketitle

\begin{abstract}
The computation of the entries of Jacobi operators associated with orthogonal polynomials has important applications in numerical analysis.  From truncating the operator to form  a Jacobi matrix, one can apply the Golub--Welsh algorithm to compute the Gaussian quadrature weights and nodes.  Furthermore, the entries of the  Jacobi operator are the coefficients in the three-term recurrence relationship for the polynomials.  This provides an efficient method for evaluating the orthogonal polynomials.  Here, we present an $\O(N)$ method to compute the first $N$ rows of Jacobi operators from the associated weight.  The method exploits the Riemann--Hilbert representation of the polynomials by solving a deformed Riemann--Hilbert problem numerically. 
  We further adapt this computational approach to  certain entire weights that are beyond the reach of current asymptotic Riemann--Hilbert techniques.  
\end{abstract}

\section{Introduction}

Infinite symmetric tridiagonal matrices with positive off-diagonal entries are referred to as Jacobi operators.  Such operators have an intimate connection to orthogonal polynomials. Given a positive weight $\dkf \rho(x) = \omega(x) \dx$ on $\mathbb R$,  $\omega$ entire with finite moments, the Gram--Schmidt procedure produces a sequence of orthogonal polynomials.  These polynomials satisfy a three-term recurrence relation, see \eqref{three-term}.  From the coefficients in the recurrence relation, a Jacobi operator
	 \begin{align}\label{jacobi-op}
J_\infty(\omega) = \begin{mat} a_0 & \sqrt{b_1} &  & && 0 \\
\sqrt{b_1} & a_1  & \sqrt{b_2} &&&\\
& \sqrt{b_2} & a_2 & \sqrt{b_3} && \\
&& \ddots & \ddots & \ddots& \\
0&&&&&
\end{mat}
\end{align} can be constructed.  Our aim is to rapidly calculate the entries of this Jacobi operator from the given weight $\omega$.



Computing the Jacobi operator \eqref{jacobi-op} is equivalent to computing the coefficients in the three-term recurrence relation \eqref{three-term}.  Gautschi states in \cite{GautschiOP}:
\begin{quote}
The three-term recurrence relation satisfied by orthogonal polynomials is arguably the single most important piece of information for the constructive and computational use of orthogonal polynomials. Apart from its obvious use in generating values of orthogonal polynomials and their derivatives, both within and without the spectrum of the measure, knowledge of the recursion coefficients (i) allows the zeros of orthogonal polynomials to be readily computed as eigenvalues of a symmetric tridiagonal matrix, and with them the all-important Gaussian quadrature rule, (ii) yields immediately the normalization coefficients needed to pass from monic orthogonal to orthonormal polynomials, (iii) opens access to polynomials of the second kind and related continued fractions, and (iv) allows an efficient evaluation of expansions in orthogonal polynomials.
\end{quote}
All applications discussed here are merely a consequence of this broad and important fact emphasized by Gautschi.

While the entries of the Jacobi operator can be calculated using the Stieljes procedure (a discretization of the  Gram--Schmidt procedure), the complexity to accurately calculate the first $N$ entries is observed to be $\O(N^3)$, see \secref{Stieljes}.  We develop an alternative to the Stieljes procedure that achieves $\O(N)$ operations: calculate the entries of the Jacobi operator using the  formulation of orthogonal polynomials as solutions to Riemann--Hilbert problems (\RHP s).   Riemann--Hilbert problems are boundary value problems for analytic functions, and over the last 20 years they have played a critical role in the advancement of asymptotic knowledge of orthogonal polynomials \cite{DeiftOrthogonalPolynomials,DeiftCollab3,DeiftCollab2,KuijlaarsRecurrence,DeiftCollab1}.  More recently, a numerical approach for solving Riemann--Hilbert problems has been developed \cite{SOPainleveII,SORHFramework}, which the authors used to  numerically calculate the $N$th and $(N-1)$st orthogonal polynomials pointwise, allowing for the computation of statistics of random matrices \cite{SOTrogdonRMT}, for varying  exponential weights of the form  $\omega(x ; N) = \E^{-N V(x)}$.  In this paper, we  extend  this numerical Riemann--Hilbert approach  to stably calculate orthogonal polynomials for non-varying weights, as well as the entries of the associated Jacobi operator.

\begin{remark}
	The Riemann--Hilbert approach 	also allows the calculation of $p_N(x)$ pointwise in $\O(1)$ operations, which can be used to achieve an $\O(N)$ algorithm for calculating Gaussian quadrature rules, {\it \'a la} the Glaser--Liu--Rokhlin algorithm \cite{GlaserLiuRokhlin,HaleTownsendQuad}, whereas the Golub--Welsch algorithm is $\O(N^2)$ operations.   While technically an $\O(N)$ algorithm, its applicability is currently limited due to large computational cost of evaluating solutions to \RHP s.  
\end{remark}

\begin{remark}
	In many cases one also considers the map from a Jacobi operator to the associated weight.  In the case that the elements of a Jacobi operator are unbounded, a self-adjoint extension $\tilde J_\infty$ is considered.  The spectral measure $\tilde \rho$ for $\tilde J_\infty$ satisfies 
\begin{align*}
 {\bf e}_1^\intercal (\tilde J_\infty-z)^{-1}{\bf e}_1 = \int_{\mathbb R} \frac{\dkf\tilde \rho(x)}{x-z}.
\end{align*}
Here ${\bf e}_1$ is the standard basis vector with a one in its first component and zeros elsewhere.  It is known  that if $\omega(x) = e^{-p(x)}$ for a  polynomial $p$ then $J_\infty$ is essentially self-adjoint (\emph{i.e.} $\tilde J_\infty$ is unique) and $\omega(x) \dx = \dkf \tilde \rho(x)$\ \cite{DeiftOrthogonalPolynomials}.  Therefore, from a theoretical standpoint, the map $\varphi$ from a subset of essentially self-adjoint Jacobi operators to this restricted class of measures is invertible.  Viewed in this way, we are concerned with the computation of $\varphi^{-1}$, the inverse spectral map.
\end{remark}

\subsection{Complexity of the Stieljes procedure}\label{sec:Stieljes}

The monic polynomials $\pi_n$ with respect to the weight $\omega$ satisfy the  three-term recurrence relation
\begin{align}\label{three-term}
\begin{split}
\pi_{-1}(x) &= 0,\\
\pi_{0}(x) &= 1,\\
\pi_{n+1}(x) &= (x-a_{n})\pi_{n}(x) - b_{n} \pi_{n-1}(x),
\end{split}
\end{align}
where
	$$a_n = \ip<x \pi_n,\pi_n> \gamma_n \qand b_n = \ip<x \pi_n,\pi_{n-1}> \gamma_{n-1} \qfor \ip<f,g> =  \int f(x) g(x) \omega(x) \dx,$$
%
where $\gamma_n$ are the normalization constants
\begin{align*}
\gamma_{n} = \ip<\pi_n,\pi_n>^{-1} = \left( \int_{\mathbb R} \pi_{n}^2(x) \omega(x)\dx \right)^{-1}.
\end{align*}
The constants $\gamma_n$ can be calculated directly from
	 $b_{n,N}$  \cite{GautschiOP}:
	$$\gamma_{n,N}^{-1}   = b_{n} b_{n-1} \cdots b_{1} \int \omega(x) \dx,$$
	though it is convenient to compute them separately from the \RHP. 

\begin{remark}
In our notation, $a_{n}$ and $b_{n}$ are the recurrence coefficients for the monic orthogonal polynomials.  The Jacobi operator that follows in  \eqref{jacobi-op} encodes the recurrence coefficients for the orthonormal polynomials.  We also remark that the orthogonal polynomials themselves can be calculated from the recurrence coefficients.  
\end{remark}


The conventional method for computing recurrence coefficients is the Stieljes procedure, which is a discretized modified Gram--Schmidt method \cite{GautschiOP}.  In this method, one replaces the inner product  associated with the weight by an $M$-point discretization:
	$$\ip<f,g>  \approx \ip<f,g>_M = \sum_{j = 1}^M w_k f(x_k) g(x_k).$$
  For example,  $\mathbb R$ can be mapped to the unit circle $\{|z| = 1\}$ using a \Mobius transformation and then the trapezoidal rule with $M$ sample points applied.   The three-term recurrence coefficients and the values of the monomial orthogonal polynomial with respect to the discretized inner product can be built up directly:
 	\meeq{
		\vc p_0^M = \vcone,\qquad \vc p_1^M = \vc x- {A_0^M  \over \ip<\vc 1,\vc 1>_M} \vc p_0^M, \qquad \vc p_{n+1}^M = \vc x \vc p_n^M - {A_{n}^M \over B_n^M} \vc p_n^M - {B_{n}^M \over B_{n-1}^M} \vc p_{n-1}^M   \ccr
		A_{n}^M = \ip<\vc x \vc p_n^M, \vc p_n^M>_M,\quad B_0^M = \ip<\vcone,\vcone>_M,\qquad B_{n}^M = \ip<\vc x \vc p_n, \vc p_{n-1}>_M,
		}
where $\vcone$ is a vector of $M$ ones, $\vc x = \vectt[x_1,\dots,x_M]$,  the product of two vectors is the pointwise product, \emph{i.e.}, $\vc x \vc p := \diag({\vc x}) \vc p$, and we use the notation $\ip<\vc f,\vc g>_M = \vc f^\intercal \diag\!\vect[w_1,\dots,w_M] \vc g$.   Noting that $B_n^M$ also equals $\ip<\vc p_n^M, \vc p_n^M>_M$,  the entries of the Jacobi operator are  approximated by
	$$a_n \approx a_n^M =  {A_n^M \over B_n^M} \qqand b_n \approx b_n^M = {B_{n}^M \over B_{n-1}^M}.$$

\begin{figure}[tbp]
\centering
\subfigure[]{\includegraphics[width=.48\linewidth]{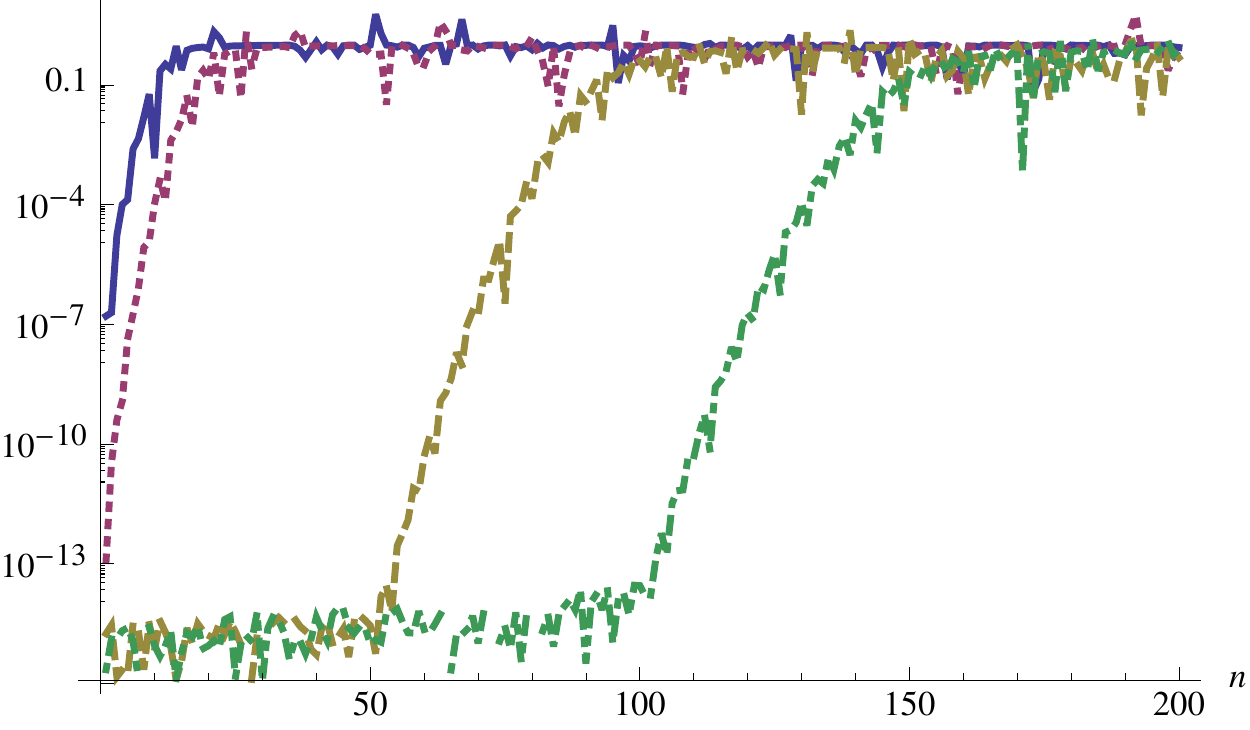}}
\subfigure[]{\includegraphics[width=.48\linewidth]{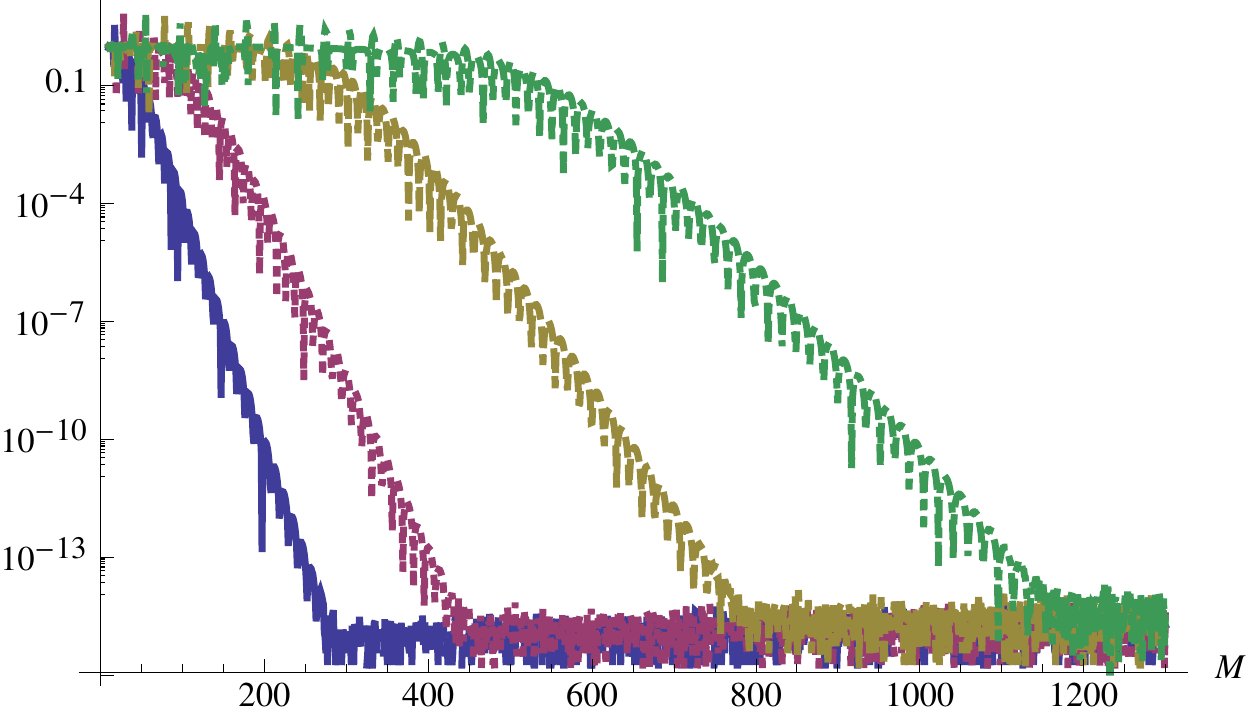}}
\caption{\label{stieljes} Relative error in approximating the Hermite polynomial recurrence coefficients. (a) Varying the number of trapezoidal rule points $M =  50$ (plain), 100 (dotted), 1000 (dashed) and 2000 (dot-dashed).  (b) The error in coefficient $n = 10$ (plain), 20 (dotted), 40 (dashed) and 60 (dot-dashed).      }
\end{figure}

  It is shown in \cite{GautschiOP} that the recurrence coefficients $a^M_{n},b^M_{n}$  obtained from this discrete inner product (via the Stieltjes procedure) converge to $a_{n},b_{n}$ for fixed $n$ as $M\rightarrow \infty$. It is stated in \cite{GautschiOP} that in $M$ must be much larger than $n$, the order of the polynomial, for $a_n^M$ and $b_n^M$ to accurately approximate $a_n$ and $b_n$.  A significant question is: how much larger must $M$ be?  In Figure~\ref{stieljes}, we show the relative error in calculating the recurrence coefficients of the Hermite weight $\E^{-x^2} \dx$ for varying $M$ and $n$ where $n$ is the index of the recurrence coefficient computed.    To accurately calculate the first 100 coefficients requires $M \approx 2000$.  Experimentally,  we observe that we require $M \sim N^2$ to achieve accuracy in calculating the first $N$ recurrence coefficients.  Each stage of the Stieljes procedure requires $\O(M)$ operations, hence the total observed complexity is $\O(N^3)$ operations.  This compares unfavorably with the $\O(N)$ operations of the \RHP\ approach.

%
%
	

\subsection{Varying weights}
	
		For stability purposes, it is necessary that we consider weights that depend on a parameter $N$ of the form 
	\begin{equation}\label{eq:Ndependent}\dkf\rho_N(x) = \omega(x ; N) \dx =  \E^{-NV_N(x)}\dx
	\end{equation}
 on $\mathbb R$, where $V_N(x) \equiv V(x)$ tends to a nice limit as $N \rightarrow \infty$.    Here $V_N(x)$ must be entire and have polynomial growth to $+\infty$ as $|x| \goto \infty$.  For a fixed $N$, we can stably calculate the $n$th row for $n \lesssim N$, hence we only consider $n = 0,1,\ldots,N$.   Varying weights of this form are useful in the study of random matrices \cite{DeiftOrthogonalPolynomials,SOTrogdonRMT}. 
 
 For Gaussian quadrature, it is more natural to consider weights of the form $\E^{-Q(x)} \dx$.   When $Q(x) = q x^m + \O(x^{m-1})$ is a degree $m$ polynomial, a simple change of variables $x \mapsto N^{1/m} x$ is sufficient to put the weight in the required form, where $V_N(x) = N^{-1} Q(N^{1/m} x) \rightarrow q x^m$ as $N \rightarrow \infty$.   
  When $Q$ is not a polynomial --- e.g. an Erd\"os weight \cite{ErdosWeight} ---  we  do a constant change of variable $x \mapsto \alpha_N x + \beta_N$,  but now $\alpha_N$ and $\beta_N$ are determined numerically in order to control the behaviour of $V_N(x) = N^{-1} Q(\alpha_N x + \beta_N)$, see Section~\ref{S:exp}.  Applying this procedure to two entire choices for $Q$, we observe a slight degradation in computational cost to an estimated count of $\O(N \log^2 N)$ operations to calculate the first $N$ rows of the Jacobi operator, see Appendix~\ref{A:Cost}.

Once we calculate the  first $N$ rows of the Jacobi operator associated with $\dkf \rho_N(x) = \E^{- Q(\alpha_N x + \beta_N)} \dx$, we can subsequently determine the first $N$ rows of   the Jacobi operator of the unscaled weight $\dkf \rho(x) = \E^{-Q(x)} \dx$.  First note that if $\pi_{j,N}$ are monic orthogonal polynomials with respect to $\dkf \rho_N$, then the monic orthogonal polynomials with respect to $\rho_N$ are $\pi_j(x) = \alpha_N^j \pi_{j,N}(\alpha_N^{-1} (x - \beta_N)).$  Thus if the entries of the Jacobi operator of $\rho_N$ are given by $a_{n,N}$ and $b_{n,N}$, then a simple calculation verifies that the entries of the original Jacobi operator of $\rho$ are
	 \begin{equation}\label{eq:recurrencetransformation}
	 	b_n = b_{n,N} \alpha_N^2\qqand a_n = a_{n,N} \alpha_N + \beta_N.
		\end{equation}
		
The factors in the right-hand sides of \eqref{eq:recurrencetransformation} depend on $N$ while the product does not.  It is advantageous to choose $N$ so that $1/c < b_{n,N} < c$ for some $c > 1$ not depending on $n$ or $N$ and the choice $n = N$ is sufficient.  The argument for this can be seen by examining the case of the Hermite polynomials where $b_n = n/2$ for $n > 0$.  Here $\alpha_{N=n}^2 = n$ so that $b_{n,N=n} = 1/2$ for all $n$ and the growth of $b_n$ is captured exactly by $\alpha_{N=n}^2$.  If $N$ is large with respect to $n$ then $b_{n,N} \ll 0$ and the computation can be prone to underflow.  Furthermore, errors accumulated in the solution of the {\RHP} are amplified by $\alpha_{N}^2$.  Our choice for $N$ is one that keeps the ``solution'' $b_{n,N}$ bounded while keeping the amplification factor $\alpha_{N}^2$ minimal.
		
\begin{remark}
	We are implicitly assuming that the equilibrium measure of the scaled potential is supported on a single interval.  This is the case when $Q$ is convex, as $V_N$ will also be convex.  
\end{remark}

%
%
%
%
%
%


\subsection{Background: Riemann--Hilbert problems}

The method presented here for the computation of recurrence coefficients relies on the numerical solution of a \emph{matrix Riemann--Hilbert problem} (\RHP). We refer the reader to \cite{SORHFramework} for a discussion of the numerical solution of {\RHP}s.  In this section we give a brief description of some components of the {\RHP} theory.  Given an oriented contour $\Gamma$, an {\RHP} consists of finding a sectionally analytic function $\Phi(z): \mathbb C \setminus \Gamma \goto \mathbb C^{2\times 2}$, depending on $n$ and $\omega(x ; N)$, such that
\begin{align*}
\lim_{\overset{z \goto x}{\text{ left of } \Gamma}} \Phi(z) = \left(\lim_{\overset{z \goto x}{z \text{ right of } \Gamma}} \Phi(z) \right) G(x ; n), ~~ G(x;n) : \Gamma \goto \mathbb C^{2\times 2}.
\end{align*}
The contour $\Gamma$ is referred to as the jump contour and $G$, as the jump matrix.  Canonically,  {\RHP}s also impose the asymptotic behaviour at $\infty$: $\lim_{|z| \goto \infty} \Phi(z)= I$.  For orthogonal polynomials, the initial {\RHP}, Problem~\ref{OP}, is reduced to canonical form via a series of transformations described in \secref{NumOPs}.


  Of course, the sense in which limits exist needs to be made precise, but this is beyond the scope of this paper, see \cite{DeiftOrthogonalPolynomials}.  We use the notation
\begin{align*}
\Phi^+(x) = \Phi_+(x) = \lim_{\overset{z \goto x}{z \text{ left of } \Gamma}} \Phi(z), ~~ \Phi^-(x) = \Phi_-(x) = \lim_{\overset{z \goto x}{z \text{ right of } \Gamma}} \Phi(z),
\end{align*}
where the $\pm$ may be located in the superscript or subscript.  Define
\begin{align*}
\CC_{\Gamma} f(z) &= \frac{1}{2 \pi i} \int_{\Gamma} \frac{f(s)}{s-z} ds,\\
\CC^\pm_{\Gamma} f(x) &= (\CC_{\Gamma} f(x))^\pm.
\end{align*}
The class {\RHP}s that we consider which satisfy $\lim_{|z| \goto \infty} \Phi(z)= I$ have the representation
\begin{align*}
\Phi(z) = I + \CC_{\Gamma}U(z),
\end{align*}
for a matrix-valued function $U \in L^2(\Gamma)$ \cite{TrogdonThesis}.  The so-called Plemelj Lemma also holds:
\begin{align*}
\CC^+_{\Gamma} f - \CC_{\Gamma}^- f = f, ~~ f \in L^2(\Gamma).
\end{align*}
From this it follows that
\begin{align*}
\Phi^+(x) = \Phi^-(x) G(x;n), ~~ \lim_{|z| \goto \infty} \Phi(z)= I,
\end{align*}
is equivalent to
\begin{align}\label{sie}
U - \CC_{\Gamma}^-U \cdot(G-I) = G-I.
\end{align}
This is a singular integral equation for $U$ and this integral equation is solved numerically with the method in \cite{SORHFramework}.


\section{Gaussian Quadrature}

In this section we discuss Gaussian quadrature on $\mathbb R$.  In Gaussian quadrature, one chooses nodes $\{x_{j}\}_{j=1}^n$ and weights $\{\omega_{j}\}_{j=1}^n$ such that
\begin{align*}
\int_{\mathbb R} p(x) \omega(x) \dx = \sum_{j=1}^n p(x_{j}) \omega_{j},
\end{align*}
for any polynomial $p(x)$ of degree $\leq 2n-1$.



\subsection{Computing nodes and weights}

We use the classical Golub--Welsh algorithm to calculate the nodes $x_j$ and weights $\omega_j$ from the Jacobi matrix, i.e.,  $n \times n$ finite section of the Jacobi operator:
\begin{align*}
J_n(\omega) = \begin{mat} a_0 & \sqrt{b_1} &  & &0 \\
\sqrt{b_1} & a_1  & \sqrt{b_2} &&\\
& \sqrt{b_2} & a_2 & \sqrt{b_3} & \\
&& \ddots & \ddots & \ddots& \\
&&&&&\sqrt{b_{n-1}}\\
0&&&& \sqrt{b_{n-1}} & a_{n-1} 
\end{mat}.
\end{align*}
The eigenvalues of $J_n(\omega)$ are the zeros of $\pi_{n}(x)$ \cite{GautschiOP}.  For each eigenvalue $\lambda_j$ of $J_n(\omega)$, an eigenvector is given by
\begin{align*}
v(\lambda_j) = (\pi_0(\lambda_j),\ldots,\pi_{n-1}(\lambda_j))^\intercal.
\end{align*}
We use $\tilde v(\lambda_j)$ to denote the normalized eigenvector, having Euclidean norm one.  Let $\mu_0 = \int \omega(x)\dx$ and it follows that $\omega_j = \mu_0 (\tilde v(\lambda_j))_1^2$, \emph{i.e.} the first component of $\tilde v(\lambda_j)$ squared times $\mu_0$ \cite{GautschiOP}. In practice, $\mu_0$ can be computed with standard quadrature routines, \emph{i.e.} Clenshaw--Curtis quadrature.


\subsection{Example: Hermite}

\begin{figure}[tbp]
\centering
\includegraphics[width=.5\linewidth]{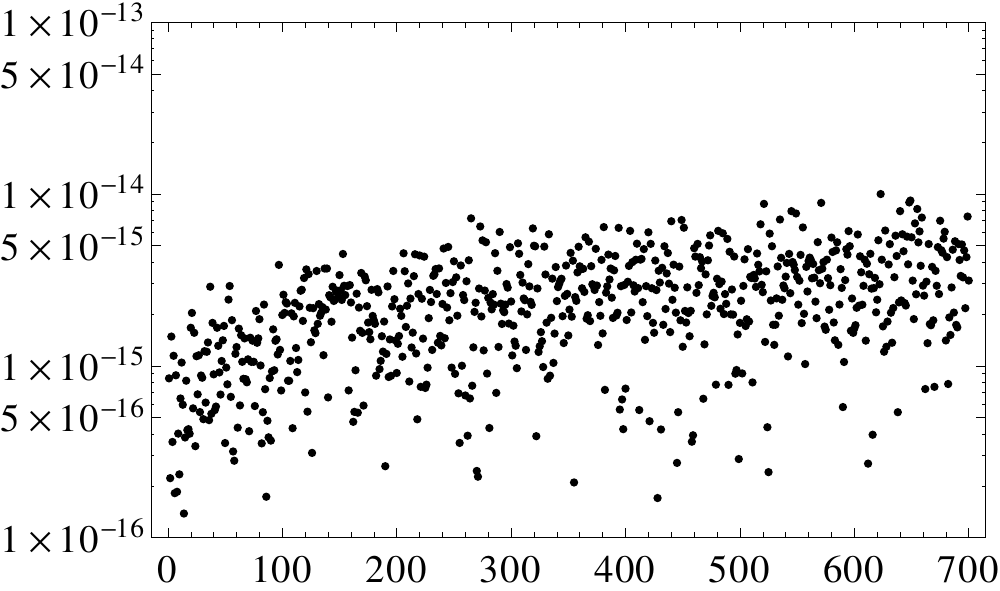}
\caption{\label{Figure:HermiteError} The relative error in $b_n$ for $n = 0,1,2,\ldots,700$.  High precision is retained.}
\end{figure}

For the Hermite polynomials, $\mu_0 = \beta_0 = \sqrt{\pi}$, $\beta_j = j/2$ for $j > 0$ and $a_j = 0$ always \cite{GautschiOP}.  We use this as a test case to demonstrate the accuracy of our computation.  In Figure~\ref{Figure:HermiteError} we show the relative error of the first 700 coefficients computed using the method described below.

\subsection{Example: $\omega(x) = e^{-x^8}$}

\begin{figure}[tbp]
\centering
\subfigure[]{\includegraphics[width=.49\linewidth]{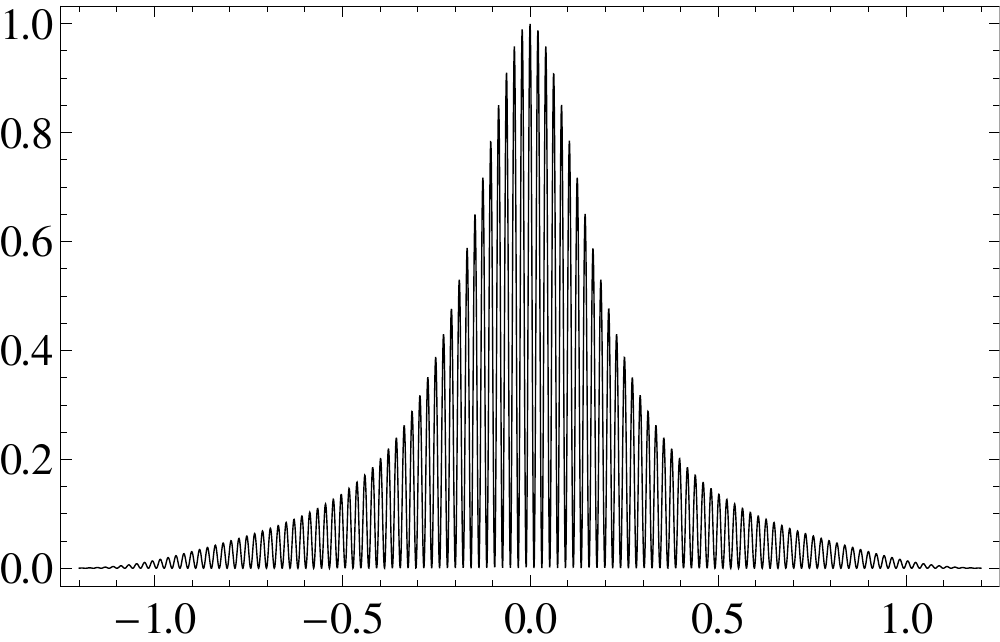}}
\subfigure[]{\includegraphics[width=.49\linewidth]{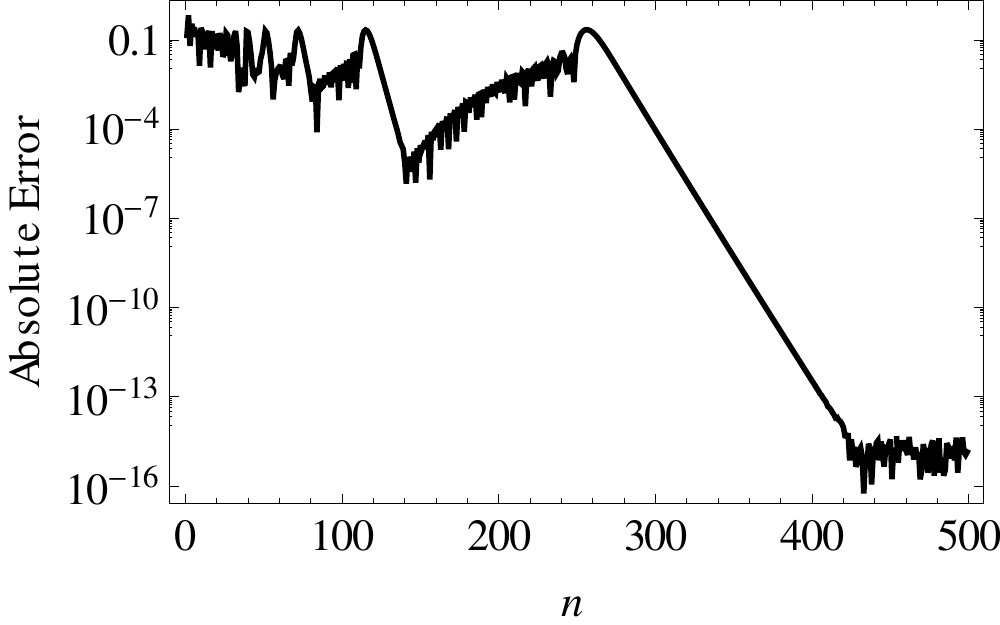}}
\caption{\label{Figure:Eight} (a) The function $f(x) \omega(x)$ for $f(x) = \cos^2(150 x)/(25x^2+1)$ and $\omega(x) = e^{-x^8}$. (b) The error in the approximation of $\int f(x) \omega(x) \dx$ with a varying number of nodes and weights.  The approximation is compared against an approximation using Clenshaw--Curtis quadrature with $10,\!000$ points.}
\end{figure}

  We demonstrate the integration of a meromorphic integrand. In all our examples, we choose to integrate positive functions so that cancellation will not increase our accuracy. See Figure~\ref{Figure:Eight} for a demonstration of this.

\subsection{Example: $\omega(x)= e^{-x^2-\sin(x)}$}

\begin{figure}[tbp]
\centering
\subfigure[]{\includegraphics[width=.49\linewidth]{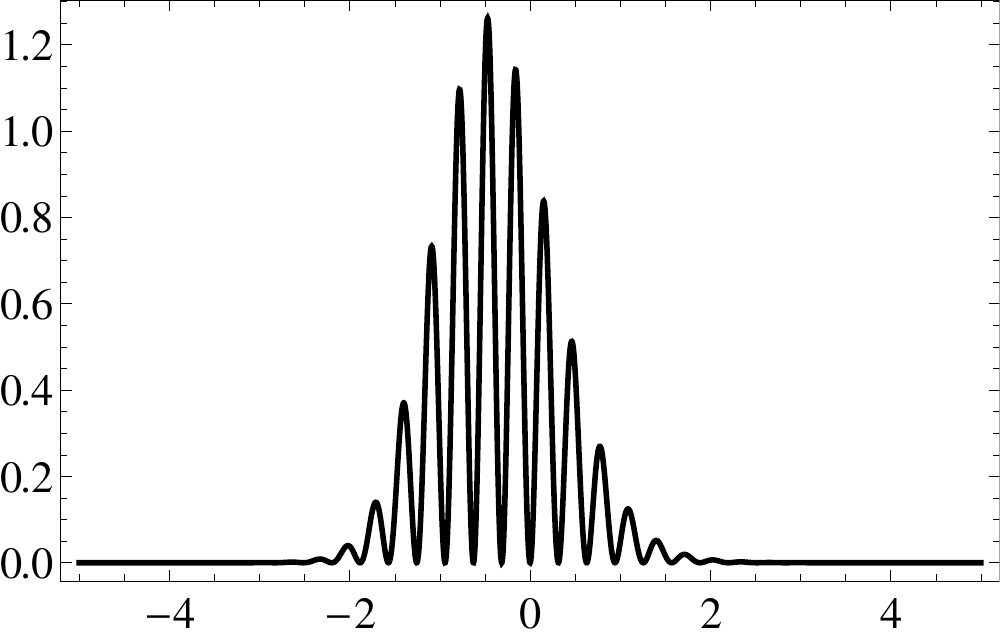}}
\subfigure[]{\includegraphics[width=.49\linewidth]{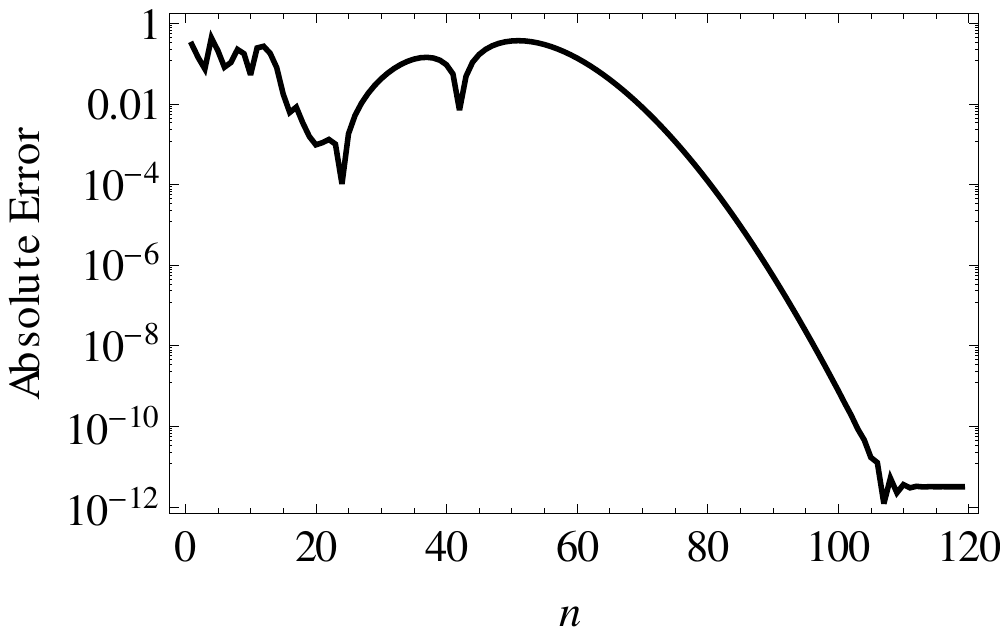}}
\caption{\label{Figure:Sheehan} (a) The function $f(x) \omega(x)$ for $f(x) = \sin^2(10x)$ and $\omega(x) = e^{-x^2-\sin(x)}$. (b) The error in the approximation of $\int f(x) \omega(x) \dx$ with a varying number of nodes and weights.  The approximation is compared against an approximation using Clenshaw--Curtis quadrature with $10,\!000$ points.}
\end{figure}

A more exotic choice is $Q(x) = x^2 + \sin(x)$.  Even though $Q$ is not a polynomial and it is not convex, numerical experiments demonstrate that the method described below to compute the recurrence coefficients is still effective.  We demonstrate using this weight with Gaussian quadrature in Figure~\ref{Figure:Sheehan}.

\begin{remark}
The only reason the method would fail on this weight would be if the equilibrium measure (see Section~\ref{sec:NumOPs}) was not supported on a single interval.
\end{remark}

\subsection{Example: $\omega(x)= e^{-\!\cosh(x)}$}

\begin{figure}[tbp]
\centering
\subfigure[]{\includegraphics[width=.49\linewidth]{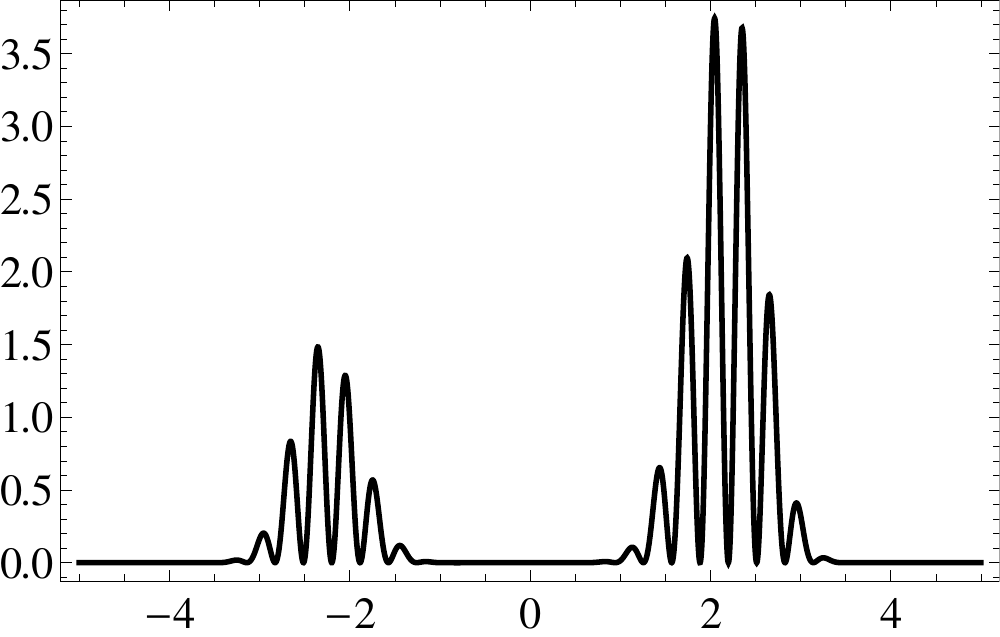}}
\subfigure[]{\includegraphics[width=.49\linewidth]{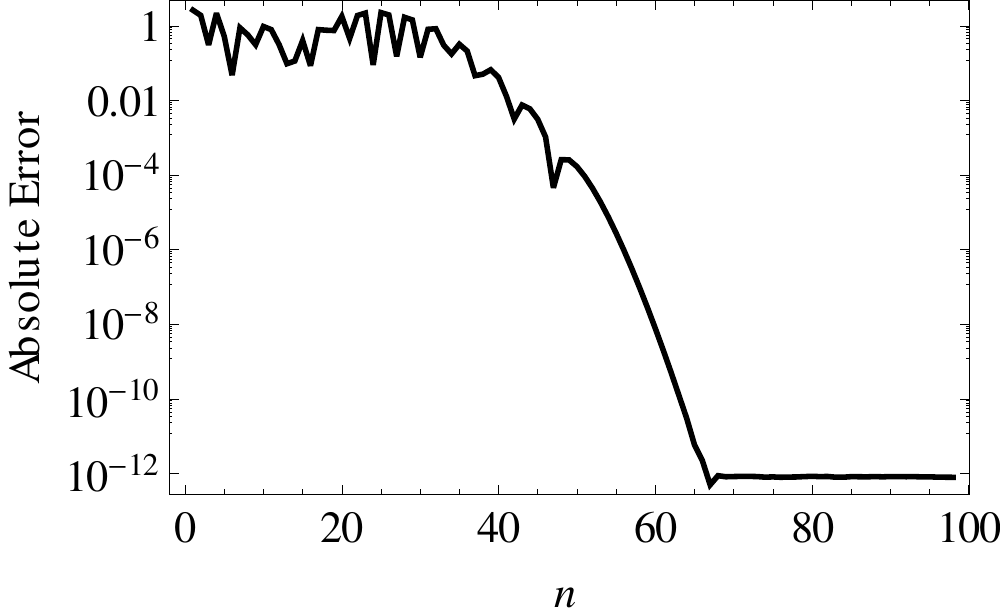}}
\caption{\label{Figure:Cosh} (a) The function $f(x) = 0.1(x^{10}+x^9)\sin^2(10x)e^{-\cosh(x)}$. (b) The error in the approximation of $\int f(x) \dx$ with a varying number of nodes and weights.  The approximation is compared against an approximation using Clenshaw--Curtis quadrature with $10,\!000$ points.}
\end{figure}

With the choice $Q(x) = \cosh(x)$, we have super-exponential decay.  We include this example to show that the method extends to $Q$ that have super-polynomial growth.  Indeed, a modification of the method is needed, see Section~\ref{S:exp}. We demonstrate using this weight with Gaussian quadrature in Figure~\ref{Figure:Cosh}.

\section{The numerical RH approach to OPs}\label{sec:NumOPs}

Here we present a concise adaptation of the numerical RH approach to computing orthogonal polynomials on $\mathbb R$ with analytic exponential weights (originally described in \cite{SOTrogdonRMT}) to calculating {\it all} orthogonal polynomials of order $n = 0,1,\ldots,N$.  For given $n$ and $N$, the \RHP\ simultaneously captures the  polynomials $\pi_{n,N}(x)$ and $\gamma_{n-1,N} \pi_{n-1,N}(x)$.  More precisely, these polynomials can be expressed in terms of the solution of an {\RHP}: 
	
\begin{problem} \label{OP} \cite{FokasOP}
The function
\begin{align*}
Y(z) =  \begin{mat}
  \pi_{n,N}(z)& \CC_{\mathbb R}{[\pi_{n,N} e^{-N V}](z)}  \\
  - { 2 \pi i} \gamma_{n-1,N} \pi_{n-1,N}(z) & - { 2 \pi i} \gamma_{n-1,N} \CC_{\mathbb R}[\pi_{n-1,N} e^{-N V}](z)
\end{mat},
\end{align*}
where
	$$\gamma_{n-1,N} = \left[ \int \pi_{n-1,N}^2(x) e^{-N V(x)} \dx\right]^{-1},$$  is the unique solution of the {\RHP} on the real line
	\begin{align*}
          Y^+(z) = Y^-(z) \begin{mat} 1 & e^{-N V(z)}\\ 0 & 1\end{mat}, \quad z \in \mathbb R, \quad Y \sim \begin{mat} z^n &0 \\ 0 & z^{-n} \end{mat}, ~~~ z \rightarrow \infty.
        \end{align*}
\end{problem}

We apply the numerical method described in \cite{SORHFramework} but conditioning steps are required. Our first task is to remove the growth/decay of $Y$ at $\infty$.  To accomplish this, we must compute the $g$-function, which   has logarithmic growth at infinity so that $e^{\pm n g(z)} \sim z^{\pm n}$, but has special jump properties so that its incorporation into the {\RHP} allows the {\RHP} to behave well in the large-$n$ limit.  

The $g$-function, which depends on $n$ and $N$, is associated with the equilibrium measure of $V(x)$:

\begin{definition}\label{EM}
	The {\it equilibrium measure} $\mu \equiv \mu_{n,N}$ is the minimizer of 
	$$\iint\log{ \frac{1}{|x - y|}} \dkf\mu(x) \dkf\mu(y) + \int \frac{N}{n} V(x) \dkf\mu(x).$$
\end{definition}
For simplicity, we assume that the equilibrium measure of $ V(x)$ is supported on a single interval $(a,b)$.  This is certainly true provided $V$ is convex \cite{DeiftOrthogonalPolynomials}.  A numerical method for computing equilibrium measures is presented in \cite{SOEquilibriumMeasure} and we leave the details to that reference. A discussion of  how the computational cost of computing the equilibrium measure varies with $N$ and $n$ is left to Appendix~\ref{A:Cost}.


From \cite{DeiftOrthogonalPolynomials} it follows that $g \equiv g_{n,N}$ is an anti-derivative of the Stieljes transform of $\mu$:
	 $$g(z) = \int^z  \int {\dkf \mu(x) \over z - x} \dz.$$
Therefore, $g$ is analytic off $(-\infty,b)$ and satisfies  
\begin{align*}
g_+(z) &+ g_-(z) = \frac{N}{n} V(z) - \ell, ~~ z \in (a,b), ~~ \ell \in \mathbb C,\\
g(z) &= \log z + \bigo(z^{-1}) ~\text{ as }~ z \goto \infty,\\
g_+(z) &- g_-(z) = 2 \pi i, ~~ z \in (-\infty,a).
\end{align*}
The method of \cite{SOEquilibriumMeasure} gives a simple representation of $g$ in terms of the Chebyshev coefficients of $V$ over $(a,b)$ and elementary functions.  Numerically,  the Chebyshev coefficients of $V'$ are replaced with a numerical approximation via the discrete cosine transform applied to $V'$ evaluated at Chebyshev points, and the $g$ associated with the resulting polynomial is calculated in closed form.  The error introduced by this approximation is uniformly small throughout the complex plane, including for $z$ close to or on $\mathbb R$.  

\subsection{Lensing the \RHP}  

We rewrite $Y(z)$ to normalize the behavior at infinity:
\begin{align}\label{Tdef}
Y(z) =  \begin{mat}
e^{n \ell \over 2} & 0\\
 0& e^{-{n \ell \over 2}}
\end{mat} T(z) \begin{mat} e^{- n g(z)} &0 \\
0& e^{n g(z)} \end{mat} \begin{mat} e^{-{n \ell \over 2}} &0 \\
0 & e^{n \ell \over 2} \end{mat},
\end{align}
so that $T(z) \sim I$ for large $z$ and $T(z)$ has a jump along the real line,  on which it satisfies
\begin{align*}
  T_+(z) &=  T_-(z) \begin{mat}
    e^{n (g_(z) - g_+(z))}  & e^{n (g_+(z) + g_-(z) + \ell  -  \frac{N}{n} V(z))}    \\
    0 & e^{n (g_+(z) - g_-(z))}\end{mat} \\
  &= T_- \begin{choices}
    \begin{mat}1 & e^{n (g_+(z) + g_-(z) + \ell  -  \frac{N}{n} V(z))} \\ 0& 1 \end{mat}   \when  z < a,\\
\\
    \begin{mat} e^{n (g_-(z) - g_+(z))}  & 1  \\  0 & e^{n (g_+(z) - g_-(z))}\end{mat} \when a < z < b,\\
\\
    \begin{mat}1 & e^{n (2 g(z) + \ell  -  \frac{N}{n}V(z))} \\ 0& 1\end{mat}  \when b < z,		
  \end{choices}
\end{align*}

\begin{figure}[htp]
\centering
\includegraphics[width=.8\linewidth]{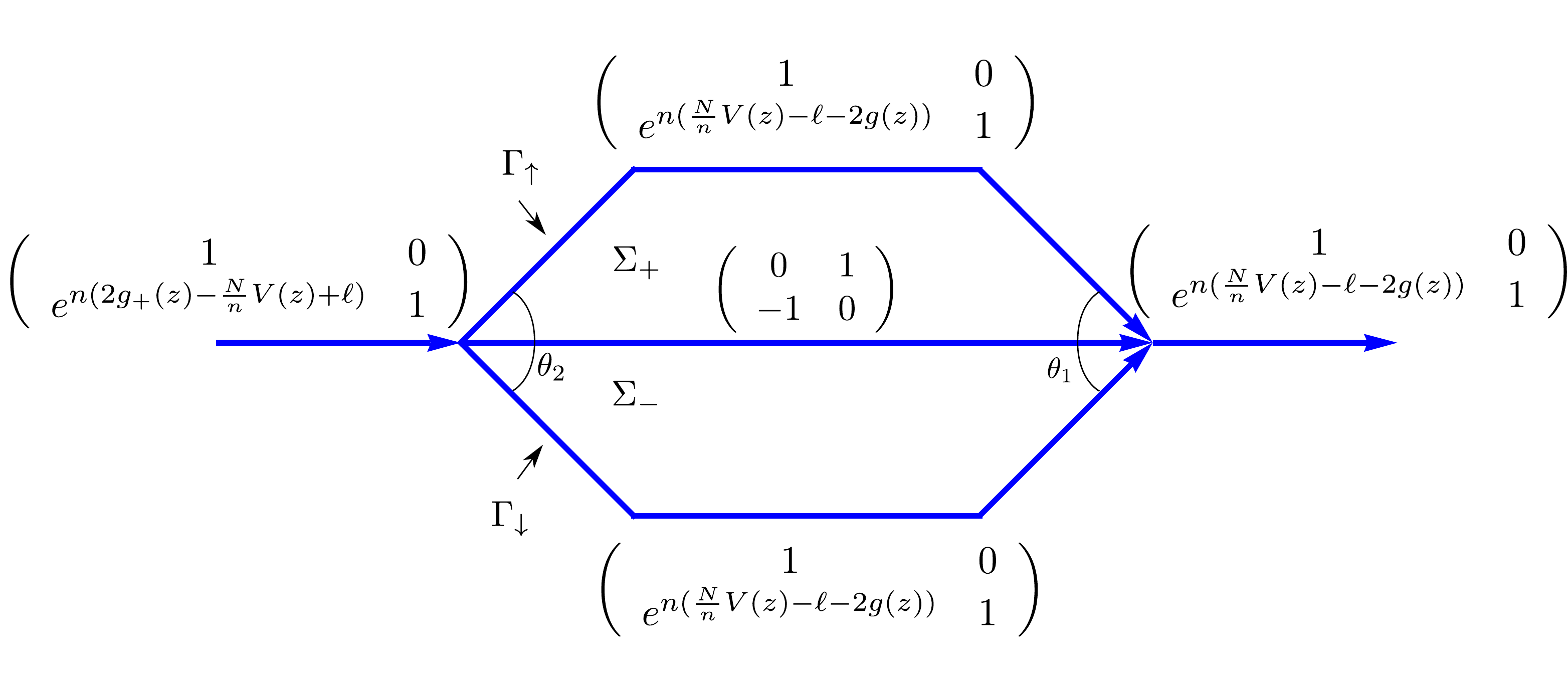}
\caption{The jumps of $S$ with the angles $\theta_1$ and $\theta_2$ labelled.\label{Figure:LensingMod}}
\end{figure}

We appeal to properties of equilibrium measures \cite[(7.57)]{DeiftOrthogonalPolynomials} to assert that 
	$$g_+(x) + g_-(x) + \ell - \frac{N}{n}V(x) < 0,$$
for $x < a$ and $x > b$, thus those contributions of the jump matrix are isolated around $a$ and $b$.  On the other hand, $g_+ - g_-$ is imaginary between $a$ and $b$ \cite[pp.~195]{DeiftOrthogonalPolynomials}, hence $e^{\pm n (g_+ - g_-)}$ is oscillatory on $(a,b)$.  The {\RHP} is deformed into the complex plane to convert oscillations into exponential decay.   A so-called lensing of the {\RHP} is introduced in Figure~\ref{Figure:LensingMod}, where we rewrite $T$ as
\begin{align}\label{Sdef}
  T(z) = S(z) \begin{choices}
    \begin{mat} 1 & 0 \\ e^{n ( \frac{N}{n} V(z) - \ell - 2 g(z))} & 1 \end{mat}, \when  z \in \Sigma_+,\\
    \begin{mat} 1 & 0 \\  e^{n ( \frac{N}{n} V(z) - \ell - 2 g(z))} & 1 \end{mat}, \when  z \in \Sigma_-,\\
    I, \otherwise .
  \end{choices}
\end{align}
By substituting
	$$g_+(z) =  \frac{N}{n} V(z) - g_-(z) - \ell,$$
within the support of $\mu$ we see that the oscillations have been removed completely from the support: 
\begin{align*}
		S_+(z) &= T_+(z) \begin{mat} 1 & 0 \\ -e^{n ( \frac{N}{n} V(z) - \ell - 2 g_+(z))} & 1 \end{mat}   = T_+(z) \begin{mat} 1 & 0 \\ -e^{n (g_-(z) -  g_+(z))} & 1 \end{mat}   \\
                &= T_-(z) \begin{mat} e^{n (g_-(z) - g_+(z))} & 1 \\ 0 & e^{n (g_+(z) - g_-(z))} \end{mat}  \begin{mat} 1 & 0 \\ -e^{n (g_-(z) -  g_+(z))} & 1 \end{mat}	\\
			 &= T_-(z) \begin{mat} 0 & 1 \\ -1  &  e^{n (g_+(z) - g_-(z))}  \end{mat}  \\
&= S_-(z) \begin{mat} 1 & 0 \\ -e^{n ( \frac{N}{n} V(z) - \ell - 2 g_-(z))} & 1 \end{mat} \begin{mat} 0  & 1 \\ -1  &  e^{n ( \frac{N}{n} V(z)-\ell -2 g_-(z))} \end{mat}  \\
			 	&= S_-(z) \begin{mat} 0& 1 \\ -1 &0 \end{mat}.
                              \end{align*}
However, we have introduced new jumps on $\Gamma_\uparrow$ and $\Gamma_\downarrow$, on which
	$$S_+(z) = T_+(z) = T_-(z) = S_-(z) \begin{mat} 1 & 0 \\ e^{n ( \frac{N}{n} V(z) - \ell - 2 g(z))} & 1 \end{mat}.$$
The choices of $\theta_1$ and $\theta_2$ in Figure~\ref{Figure:LensingMod} are discussed below.

\subsection{Removing the connecting jump}  \label{sec:Remove}

The above process has shifted oscillations to exponential decay.  However, the goal of these deformations of the {\RHP} is  to maintain accuracy of the numerical algorithm for large $n$.  Following the heuristics set out in \cite{SOTrogdonNNSD} we must isolate the jumps to neighborhoods of the endpoints $a$ and $b$ by removing the jump along $(a,b)$.  We introduce a {\it global parametrix} that solves the {\RHP} exactly on this single contour.  In other words, we require  a function which satisfies the following {\RHP}:
	$$\mathcal N_+(x) = \mathcal N_-(x) \begin{mat} 0 & 1 \\ -1 & 0 \end{mat}, ~\text{ for }~ a < x < b ~\text{ and }~ \mathcal N(\infty) = I.$$
The solution is \cite{DeiftOrthogonalPolynomials}
	$$\mathcal N(z)  = {1 \over 2 \nu(z) } \begin{mat} 1 &  i \\
	- i & 1 \end{mat} + {\nu(z) \over 2}  \begin{mat} 1 &  - i \\
	 i & 1  \end{mat} 
	~\text{ for } \nu(z) = \left({z - b \over z - a}\right)^{1/4};$$
\emph{i.e.}, $\nu(z)$ is a solution to
	$$\nu_+(x) = i \nu_-(x) ~\text{ for }~ a < x < b ~\text{ and }~ \nu(\infty) = 1.$$

An issue with using $\mathcal N(z)$ as a parametrix is that it introduces singularities at $a$ and $b$, hence we also  introduce {\it local parametrices} to avoid these singularities.  Generically, the equilibrium measure $\psi(x)$ has {\it square-root decay} at the endpoints of its support, and so parametrices based on Airy functions could be employed \cite{DeiftOrthogonalPolynomials}.  However, it is possible that the equilibrium measure has higher-order decay at the endpoints of its support, in which case the local parametrices are given only in terms of a {\RHP} \cite{HigherOrderTracyWidom}. To avoid this issue, we introduce the {\it trivially constructed local parametrices} which satisfy the jumps of $S$ in neighborhoods of $a$ and $b$:
\begin{align*}
P_a(z) = \left( \begin{choices}
	\begin{mat}1 & 0 \\ 1 & 1\end{mat} \when {\pi -\theta_2} < \arg(z-a) < \pi \\
	\begin{mat}1& -1 \\ 1 & 0\end{mat} \when -\pi < \arg(z-a) <{ -\pi+\theta_2}  \\
	\begin{mat}0& -1 \\ 1 & 0\end{mat} \when {-\pi + \theta_2 }< \arg(z-a) <0  \\	
	I \otherwise
	\end{choices}\right) \begin{mat} e^{n ( \frac{N}{n} V(z) - \ell - 2 g(z))} & 0 \\ 0 & e^{-n ( \frac{N}{n} V(z) - \ell - 2 g(z))} \end{mat},
\end{align*}
and
\begin{align*}
P_b(z) = \left(\begin{choices}
	\begin{mat}1 & 0 \\ -1 & 1\end{mat} \when {\theta_1} < \arg(z-b) < \pi \\
	\begin{mat}0 & -1 \\ 1 & 1\end{mat} \when -\pi < \arg(z-b) <-{ \theta_1}  \\
	\begin{mat}1 & -1 \\ 0 & 1\end{mat} \when -{ \theta_1} < \arg(z-b) < 0   \\	
	I \otherwise
	\end{choices} \right)\begin{mat} e^{n ( \frac{N}{n} V(z) - \ell - 2 g(z))} & 0 \\ 0& e^{-n ( \frac{N}{n} V(z) - \ell - 2 g(z))}\end{mat},
\end{align*}
where the angles $\theta_1$ and $\theta_2$ are determined in \secref{contourscaling}. 

\begin{figure}[htp]
\centering
\includegraphics[width=\linewidth]{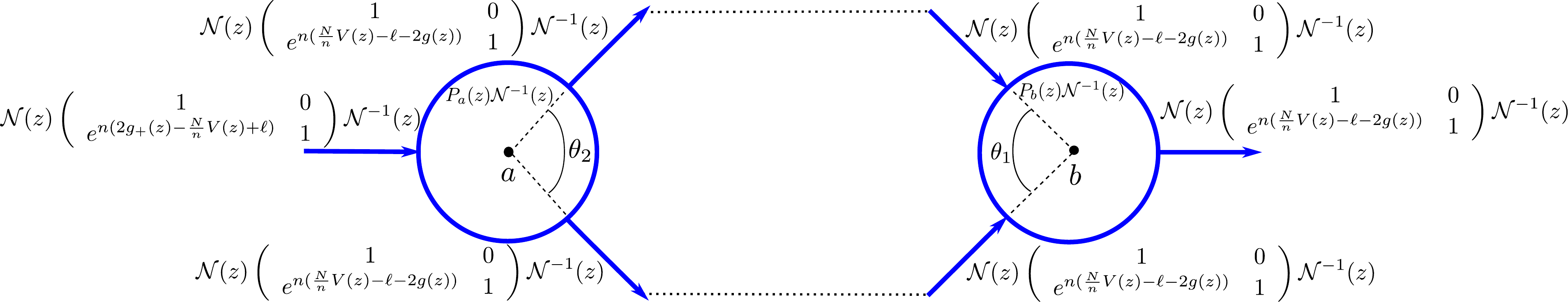}
\caption{The jumps of $\Phi$. We use a counter-clockwise orientation on the circles about $a$ and $b$. \label{Figure:PhiRHP}}
\end{figure}

We write
\begin{align}\label{Phidef}
  S(z) = \Phi(z) \begin{choices} \mathcal N(z)  \when |z - a| > r \hbox{ and } |z - b| > r, \\
    P_b(z) \when |z - b| < r, \\
    P_a(z) \when |z - a| < r. \end{choices}
\end{align}
The final {\RHP} for $\Phi$ satisfies the jumps depicted in Figure~\ref{Figure:PhiRHP}.  Note that in general $r$, the radius of the circles in Figure~\ref{Figure:PhiRHP}, depends on $N$.  We discuss this in more detail in the next section.

In practice, we do not use infinite contours.  We truncate contours when the jump matrix is to machine precision the identity matrix and this can be justified following arguments in \cite{SOTrogdonNNSD}.  In all cases we consider here, after proper deformations the jump matrices are infinitely smooth  and are exponentially decaying to the identity matrix for large $z$.

\subsection{Contour scaling and truncation}\label{sec:contourscaling}

\begin{figure}[htp]
\centering
\includegraphics[width=.9\linewidth]{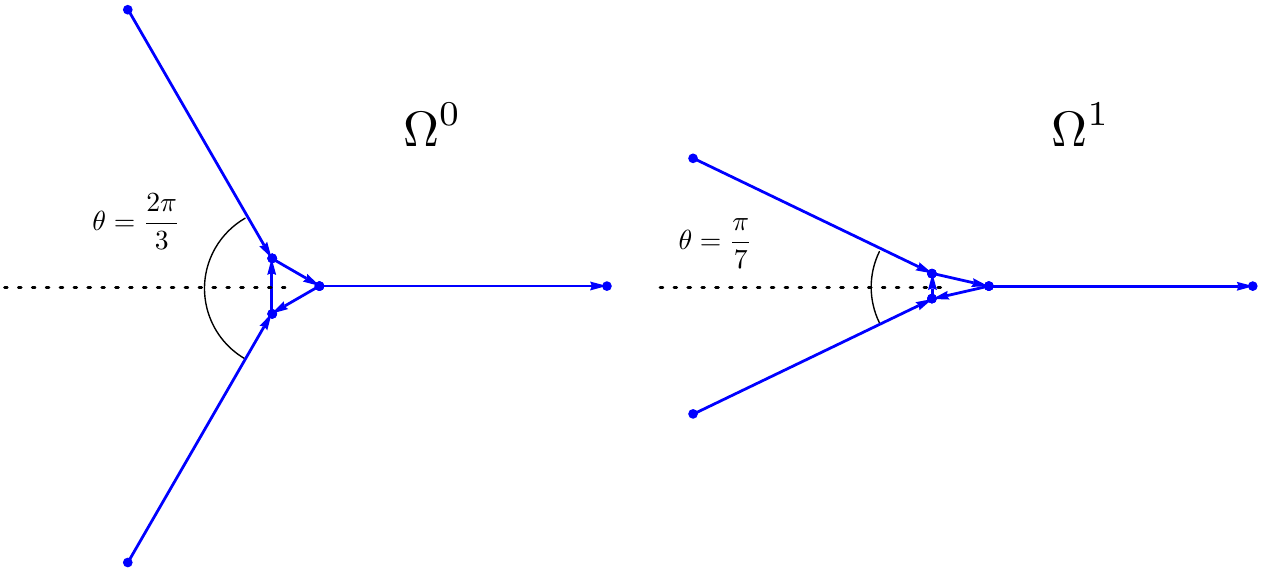}
\caption{The pre-scaled $\Omega^0$ used for non-degenerate endpoints and the pre-scaled $\Omega^1$ used for first-order degenerate endpoints.  These are the contours that are used in practice. \label{Figure:DegenerateOmega}}
\end{figure}

Evaluating the functions that appear in the jump matrix of the {\RHP} for $\Phi$ is not expensive due to the efficiency of the method for computing $g(z)$.  This allows us to perform an adaptive truncation of the jump matrices.  Before we discuss this, we must discuss the radius of the circles in Figure~\ref{Figure:PhiRHP}.

 In the case of a non-degenerate (generic case) equilibrium measure, $ \frac{N}{n} V(z)- \ell -2g(z) \sim c_a (z-a)^{3/2}$ as $z \goto a$ and $ \frac{N}{n} V(z) - \ell -2g \sim c_b (z-b)^{3/2}$ as $z \goto b$.  In this case,  we scale the contours like $n^{-2/3}$:
 	$$\Omega_n^1 = -n^{-2/3} \Omega^0 + a ~\text{ and }~ \Omega_n^2 = n^{-2/3} \Omega^0 + b,$$
for a fixed domain $\Omega^0$ that is used in practice is depicted in the left graph of Figure~\ref{Figure:DegenerateOmega}, and the angle of the contours are chosen to match the direction of steepest descent.  This implies that $r \sim n^{-2/3}$.
 In the first-order degenerate case (e.g., $V(x) = {x^2}/5 - {4} x^3/15 + {x^4}/20 + {8x}/5$), $ \frac{N}{n} V(z)- \ell -2g(z) \sim c_b (z-b)^{7/2}$ as $z \goto b$ and we scale like $n^{-2/7}$ (implying $r \sim n^{-2/7}$) at the degenerate endpoint:    
 	$$\Omega_n^1 = -n^{-2/3} \Omega^0 + a ~\text{ and }~ \Omega_n^2 = n^{-2/7} \Omega^1 + b,$$
where $\Omega^1$ is depicted in the right graph of Figure~\ref{Figure:DegenerateOmega} (the angle is sharper to attach to the new direction of steepest descent).  In general, $\theta_1$ and $\theta_2$ in Figure~\ref{Figure:LensingMod} are chosen based on the angles for $\Omega^0$ and $\Omega^1$.  Higher order degenerate points can only take the form  $ \frac{N}{n} V(z)- \ell -2g(z) \sim c_b (z-b)^{(3 + 4 \lambda)/2}$ for integer $\lambda$ \cite{HigherOrderTracyWidom}, and can be treated similarly by scaling a suitably constructed contour like $n^{-2/(3 + 4 \lambda)}$.

Each non-self-intersecting component $\Omega$ of the contour in Figure~\ref{Figure:PhiRHP} is given by a transformation $\Omega = M([-1,1])$ \cite{SORHFramework}.  After truncation, it is often convenient to replace arcs with straight lines.  Therefore, we assume that each $M(x) = \alpha x + \beta$ is just an affine transformation.  Given a jump matrix $G$ defined on $\Omega$  we wish to truncate the contour when $\|G(k)-I\| < \epsilon$.  

\begin{algorithm}\label{truncate}
\textbf{Contour truncation}
\begin{enumerate}
\item Fix $N_{\text{\rm grid}} > 0$  and consider the values
\begin{align*}
\|G(k_1)-I\|,\|G(k_2)-I\|,\ldots,\|G(k_{M_{\text{\rm grid}}})-I\|,\\
k_j = M\left( 2 \frac{j}{M_{\text{\rm grid}}} -1 \right),\quad M_{\text{\rm grid}}= \lceil (|\alpha|+1) N_{\text{\rm grid}} \rceil.
\end{align*}
Here $\{k_j\}$ represents a uniform grid on $\Omega$ that increases in number of points as $\Omega$ increases in extent.  There is a simple ordering on $\Omega$: we say $s < t$, $s,t \in \Omega$ if and only if $M^{-1}(s) < M^{-1}(t)$. 

\item If $\|G(k_1)-I\| \geq \epsilon$ set $S = \{k_1\}$ otherwise $S =\{\}$.
\item For $j=1,\ldots,M_{\text{\rm grid}}-1$, if $(\|G(k_j)-I\|-\epsilon)(\|G(k_{j+1})-I\|- \epsilon) < 0$, indicating a sign change, add $(k_j+k_{j+1})/2$ to the set $S$.
\item If $S$ has an odd number of elements add $k_{M_{\text{\rm grid}}}$ to $S$.
\item Order $S = \{t_1,\ldots,t_{2n}\}$ according to $t_j < t_{j+1}$.  Replace $\Omega$ with the contours $\{\Omega_1,\ldots,\Omega_n\}$ where $\Omega_j$ is a straight line that connects $t_{2j-1}$ and $t_{2j}$.
\end{enumerate}
\end{algorithm}

\begin{figure}[tbp]
\centering
\subfigure[]{\includegraphics[width=.5\linewidth]{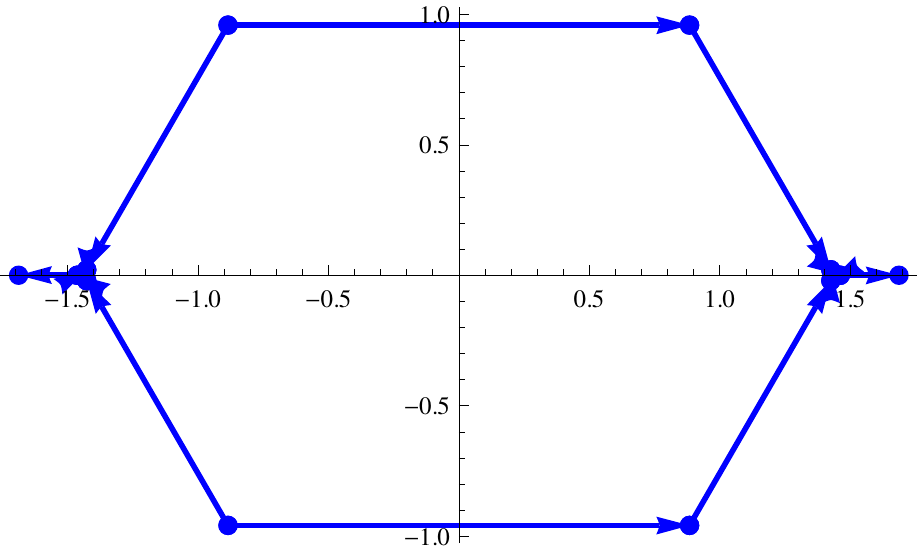}}
\subfigure[]{\includegraphics[width=.5\linewidth]{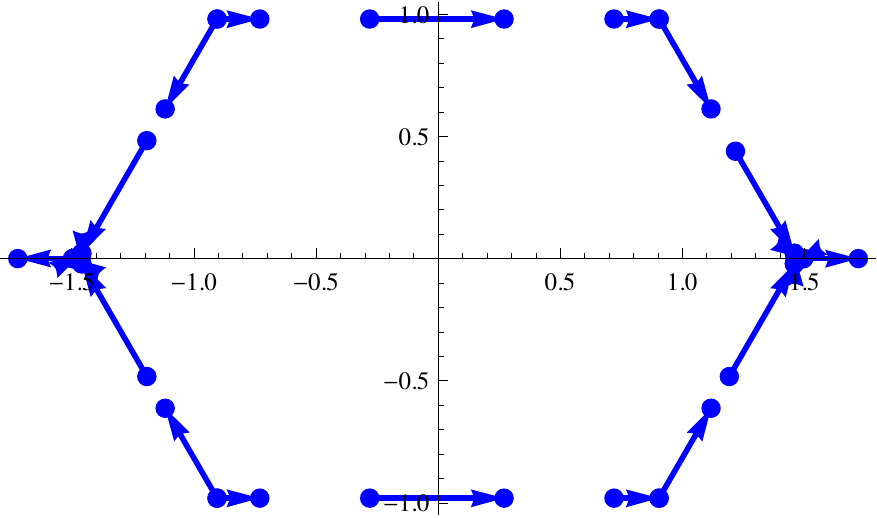}}
\subfigure[]{\includegraphics[width=.5\linewidth]{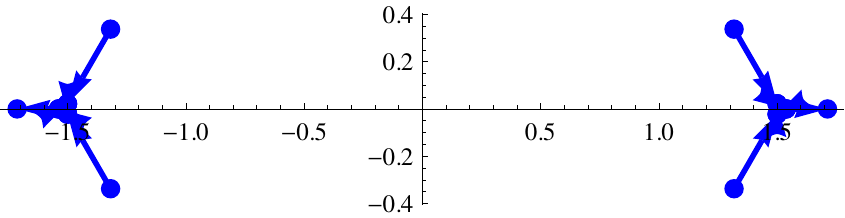}}
\caption{A demonstration of the resulting contours after the application of Algorithm~\ref{truncate} with $V(x) = x^8$ and $N =1$: (a) $n = 20$, (b) $n=24$, (c) $n=30$. \label{Figure:truncate}}
\end{figure}

Provided that $N_{\text{\rm grid}}$ is sufficiently large so that at most one zero of $\|G(k_j)-I\|-\epsilon$ lies between successive points $k_j$ and $k_{j+1}$ the algorithm produces a consistent truncation of the contour in the sense that 
\begin{align*}
\|G(k)-I\| \leq \sqrt{\frac{\sup \partial_k \|G(k)-I\|^2}{(|\alpha|+1) N_{\rm grid}} + \epsilon^2} \text{  for } k \in \Omega \setminus \{\Omega_1 \cup \cdots \cup \Omega_n\}.
\end{align*}
If $\|G(k)-I\|$ happens to oscillate wildly as $n$ increases then $N_{\text{\rm grid}}$ may have to increase with respect to $n$.  In practice, $\|G(k)-I\|$ is seen to either be monotonic or vary slowly.  Indeed, this is the purpose of the deformation of an {\RHP}.  In both cases, the algorithm produces a good truncation.  As a product of this, contours are generally divided up into at most three components after truncation, see Figure~\ref{Figure:truncate}.  The number of mapped Chebyshev points on $\Omega$ used to discretize the singular integral equation is used on each of the subdivided contours to ensure enough points. In practice, this procedure ensures accuracy for all $n$ and $N$.


\subsection{Numerically solving the {\RHP}}

The {\RHP} for $\Phi$, after scaling and truncation is solved numerically.  To avoid too much complication, we only describe the essence of the method \cite{SORHFramework} that is applied. This method discretizes \eqref{sie} via collocation using an appropriately defined Chebyshev basis.  As is described in \cite{SOTrogdonRMT} the scaling and truncation is sufficient to ensure that the number of collocation points needed to compute the solution accurately is bounded (with respect to $n$).  Thus to solve the {\RHP} for $\Phi$ the computational cost is $\bigo(1)$.  Once $\Phi$ is calculated, we can evaluate orthogonal polynomials pointwise by undoing the sequence of transformations to determine $Y$.  However, we will see that the recurrence coefficients are also calculable directly.

\subsection{Computing recurrence coefficients}

For a given  $V(x)$, we discuss how to accurately obtain $a_{n,N}$ and $b_{n,N}$.  For large $z$, the solution $\Phi$ of the {\RHP} in Figure~\ref{Figure:PhiRHP} is given by
\begin{align*}
\Phi(z) = T(z) \mathcal N^{-1}(z) = \begin{mat} e^{-n\ell/2} & 0 \\ 0 & e^{n\ell/2} \end{mat} Y(z) \begin{mat} e^{n\ell/2} & 0 \\ 0 & e^{-n\ell/2} \end{mat} \begin{mat} e^{n g(z)} & 0 \\ 0 & e^{-n g(z)} \end{mat} \mathcal N^{-1}(z).
\end{align*}
For fixed $N$, we wish to compute the polynomials, orthogonal with respect to $e^{-NV(x)}\dx$ of order $n=0,1,\ldots,N$.  From \cite{KuijlaarsRecurrence} we know that
\begin{align*}
b_{n,N} = (Y_{1}(n))_{12} (Y_1(n))_{21}, ~~ a_{n,N} = \frac{(Y_2(n))_{12}}{(Y_1(n))_{12}} - (Y_1(n))_{22},\\
Y(z) = \left( I + \frac{1}{z} Y_1(n) + \frac{1}{z^2} Y_2(n)  + \bigo \left( \frac{1}{z^3} \right) \right) \begin{mat} z^{n} & 0 \\ 0 & z^{-n} \end{mat}.
\end{align*}
It is inconvenient to compute $Y_2(n)$ so we look for alternate formulae for $a_{n,N}$ and $b_{n,N}$.  From \cite[p.~47]{DeiftOrthogonalPolynomials}
\begin{align*}
a_{n,N} &= (Y_1(n))_{11}-(Y_1(n+1))_{11},\\
b_{n,N} &= (Y_1(n))_{12}(Y_1(n))_{21}.
\end{align*}
Computing $Y_1(n+1)$ it is more accurate than computing $Y_2(n)$.  It is also more efficient when the sequence $a_{j,N}$, $j=0,1,2,\ldots$ is computed with the computation reused at each step.


To compute $Y_1(n)$ we use the formula
\begin{align*}
I + \frac{Y_1(n)}{z} &+ \bigo\left(\frac{1}{z^2}\right) \\
&= \begin{mat} e^{n\ell/2} & 0 \\ 0 & e^{-n\ell/2} \end{mat} \left( I + \frac{\Phi_1(n)}{z} + \bigo\left(\frac{1}{z^2}\right) \right) \left( I + \frac{\mathcal N_1}{z}+ \bigo\left(\frac{1}{z^2}\right) \right)\\
&\times\left( I +  \begin{mat} - n g_1 & 0 \\ 0 & n g_1 \end{mat} \frac{1}{z} + \bigo\left(\frac{1}{z^2}\right) \right) \begin{mat} e^{-n\ell/2} & 0 \\ 0 & e^{n\ell/2} \end{mat},\\
\Phi(z) &= I +  \frac{\Phi_1(n)}{z} + \bigo\left(\frac{1}{z^2} \right),\\
\mathcal N(z) &= I +  \frac{\mathcal N_1(n)}{z} + \bigo\left(\frac{1}{z^2} \right),\\
g(z) &=\log z + \frac{g_1}{z} + \bigo\left( \frac{1}{z^2} \right).
\end{align*}
Therefore
\begin{align*}
Y_1(n) = \begin{mat} e^{n\ell/2} & 0 \\ 0 & e^{-n\ell/2} \end{mat} \left( \Phi_1(n) + \mathcal N_1 + \begin{mat} - n g_1 & 0 \\ 0 & n g_1 \end{mat}\right) \begin{mat} e^{-n\ell/2} & 0 \\ 0 & e^{n\ell/2} \end{mat}
\end{align*}

Here $\mathcal N_1$ can be computed explicitly in terms of the endpoints of the support of the equilibrium measure. $\Phi_1$ is determinable via a contour integral once the {\RHP} for $\Phi$ is solved.  Recall, the method in \cite{SORHFramework} returns $U$ such that 
\begin{align*}
I + \frac{1}{2 \pi i} \int_{\Gamma} \frac{U(s)}{s-z} ds \approx \Phi(z).
\end{align*}
Then
\begin{align*}
\Phi_1(n) \approx - \frac{1}{2 \pi i} \int_{\Gamma} U(s) ds.
\end{align*}
Finally, 
\begin{align*}
 g_1 = \frac{1}{2 \pi i} \int_{a}^b x \dkf\mu(x).
\end{align*}

When $V(x)$ is polynomial the computational cost to compute any given recurrence coefficient $\bigo(1)$.  When $V(x)$ is entire, but non-polynomial, we obtain conjectures on the computational cost.  A full discussion of this is found in Appendix~\ref{A:Cost}.

\subsection{Regularizing exponential growth}\label{S:exp}

The method in the previous section fails for non-polynomial $Q$, for example $Q(x) = \cosh(x)$.  For polynomial $Q$, the effect of $Q(x) \mapsto V_N(x)$ is to force the support of the equilibrium measure to be $\bigo(1)$.  Thus, when $Q(x)$ grows faster than any polynomial it is no longer clear a priori what the correct scaling is.  We look for functions $\alpha_N> 0$ and $\beta_N$ such that the equilibrium measure for 
\begin{align*}
V_N(x) = N^{-1} Q(\alpha_N x + \beta_N ),~~ N > 1,
\end{align*}
has support $[-1,1]$.  A simple root-finding method can be set up to find $\alpha_N$ and $\beta_N$, and if $Q$ symmetric then $\beta_N = 0$.  This indicates that the choice $m = \log(N)/\log(\alpha_N)$ ($m =1$ when $N = 1$) in the previous section is correct choice for exponentially growing potentials.

\begin{remark}
  The resulting method is very stable to small perturbations of $\alpha_N$ and $\beta_N$. Indeed, any value of $\alpha_N$ or $\beta_N$ such that the support of the equilibrium measure is $\bigo(1)$ is an acceptable scaling. 
\end{remark}

\begin{remark}
The constants $\alpha_N$ and $\beta_N$ provide the asymptotic behavior of the coefficients of the Jacobi operator for large $n$ of the unscaled Jacobi operator  without requiring the solution of any {\RHP}, see \cite[Theorem 15.2]{LevinLubinskyOPs}. 


\end{remark}

\section{Conclusion}

We have demonstrated a new, stable and accurate algorithm for computing the entries of Jacobi operators associated with polynomials orthogonal on $\mathbb R$ with respect to analytic, exponential weights. Importantly, the computational cost to compute $n$ rows of the Jacobi operator is $\bigo(N)$ when $V(x)$ is a polynomial.  We estimate the computational cost being $\bigo(N \log^2 N)$ when $V(x)$ is entire for our examples.  As a trivial consequence computing the coefficients in the classical three-term recurrence relation for these polynomials allows us to compute the polynomials themselves.  Applying the Golub--Welsh algorithm, we compute Gaussian quadrature weights and nodes accurately.

	This approach can be extended to measures supported on a single interval, using the {\RHP} derived in \cite{kuijlaars2004riemann}.  This setting is simplified for fixed weights as the equilibrium measure does not depend on the weight: it is always an arcsine distribution.  Another extension is to weights whose equilibrium measure is supported on multiple intervals, see \cite{DeiftCollab1}.  The numerical approach for equilibrium measures is applicable to the multiple interval setting \cite{SOEquilibriumMeasure}, however, the Newton iteration underlying the determination of the support can have difficulties at transition points, where the number of intervals of support grows or shrinks.   
	

\appendix

\section{Computational cost}\label{A:Cost}

The arguments presented in \cite{SOTrogdonRMT} demonstrate that once $g(z)$ is computed, the computational cost required to construct and solve the {\RHP} is bounded.  Therefore, it is necessary to understand the computational cost associated with computing the equilibrium measure for varying $N$ and $n$ and when $V$ depends on $N$.  The method in \cite{SOEquilibriumMeasure} relies on two approximations.  The first is the representation of $\frac{N}{n} V'(x)$ in a Chebyshev series.  The second is the convergence of Newton's method to determine the support of the equilibrium measure. Thus computational cost behaves like $\bigo(K k + k \log k)$ where $K$ is the number of Newton iterations and $k$ is the number of coefficients in the Chebyshev series that are needed to approximate $V'(x)$ to machine precision.

In practice, we see that $K = \bigo(1)$.  By the chosen scalings, the endpoints of the equilibrium measure are $\bigo(1)$ and an initial guess of $[-1,1]$ is sufficient.  Additionally, one could use an initial guess from a different $N$ or $n$ value as an initial guess to provide a minor improvement in the computation time.

When $V$ is itself an $m$th order polynomial then it follows that $k = m-1$ is bounded and the computational cost to compute $g(z)$, construct the {\RHP} and solve the {\RHP} is truly $\bigo(1)$.  This issue is much more delicate when considering, say, $V_N(x) = x^2 + 1/N \sin(\sqrt{N} x)$ (which corresponds to a normalized $V(x) = x^2 + \sin(x)$) or $V_N(x) = N^{-1} \cosh(d_N x)$ (which corresponds to a normalized $V(x) = \cosh(x)$).  In the first case, $V_N'''(x)$ is unbounded which forces $k$ to increase with respect to $N$.  In both of these examples, we observe experimentally that $k = \bigo(\log N)$, implying that the  computational cost for $V(x) = x^2 + \sin(x)$ is $\bigo(N \log^2N)$.  It is also observed that it takes a bounded amount of time to compute $d_N$ and a bounded number of Chebyshev coefficients to compute the equilibrium measure for $V_N(x) = N^{-1} \cosh( d_N x)$ and therefore it is conjectured that the full computation (computing $g(z)$, constructing the {\RHP} and solving the {\RHP}) is $\bigo(1)$  in this case.  Further details are  beyond the scope of this paper.

\bibliographystyle{plain}
\bibliography{GaussianQuadrature}

\begin{thebibliography}{10}

\bibitem{HigherOrderTracyWidom}
T.~Claeys, I.~Krasovsky, and A.~Its.
\newblock Higher-order analogues of the {Tracy--Widom} distribution and the
  {Painlev\'e} {II} hierarchy.
\newblock {\em Comm. Pure Appl. Math.}, 63:362--412, 2010.

\bibitem{DeiftOrthogonalPolynomials}
P.~Deift.
\newblock {\em Orthogonal Polynomials and Random Matrices: a Riemann--Hilbert
  Approach}.
\newblock AMS, 2000.

\bibitem{DeiftCollab3}
P.~Deift, T.~Kriecherbauer, K.~T.-R. McLaughlin, S.~Venakides, and X.~Zhou.
\newblock Asymptotics for polynomials orthogonal with respect to varying
  exponential weights.
\newblock {\em Internat. Math. Res. Notices}, 16:759--782, 1997.

\bibitem{DeiftCollab2}
P.~Deift, T.~Kriecherbauer, K.~T.-R. McLaughlin, S.~Venakides, and X.~Zhou.
\newblock Strong asymptotics of orthogonal polynomials with respect to
  exponential weights.
\newblock {\em Comm. Pure Appl. Math.}, 52(12):1491--1552, 1999.

\bibitem{DeiftCollab1}
P.~Deift, T.~Kriecherbauer, K.~T.-R. McLaughlin, S.~Venakides, and X.~Zhou.
\newblock Uniform asymptotics for polynomials orthogonal with respect to
  varying exponential weights and applications to universality questions in
  random matrix theory.
\newblock {\em Comm. Pure Appl. Math.}, 52(11):1335--1425, 1999.

\bibitem{ErdosWeight}
P.~Erd\"os.
\newblock On the distribution of roots of orthogonal polynomials.
\newblock In G.~Alexits et~al., editor, {\em Proceedings of the Conference on
  the Constructive Theory of Functions}, pages 145--150. Akademiai Kiado,
  Budapest, 1972.

\bibitem{FokasOP}
A.~S. Fokas, A.~R. Its, and A.~V. Kitaev.
\newblock The isomonodromy approach to matric models in 2d quantum gravity.
\newblock {\em Comm. Math. Phys.}, 147(2):395--430, 1992.

\bibitem{GautschiOP}
W.~Gautschi.
\newblock {\em Orthogonal Polynomials: Applications and Computation}.
\newblock Oxford University Press, 2004.

\bibitem{GlaserLiuRokhlin}
A.~Glaser, X.~Liu, and V.~Rokhlin.
\newblock A fast algorithm for the calculation of the roots of special
  functions.
\newblock {\em SIAM J. Sci. Comp.}, 29(4):1420--1438, 2007.

\bibitem{HaleTownsendQuad}
N.~Hale and A.~Townsend.
\newblock Fast and accurate computation of {G}auss--{L}egendre and
  {G}auss--{J}acobi quadrature nodes and weights.
\newblock {\em SIAM J. Sci. Comp.}, 35(2):A652--A674, 2013.

\bibitem{KuijlaarsRecurrence}
A.~B.~J. Kuijlaars and P.~M.~J. Tibboel.
\newblock The asymptotic behaviour of recurrence coefficients for orthogonal
  polynomials with varying exponential weights.
\newblock {\em J. Comput. Appl. Math.}, 233(3):775--785, 2009.

\bibitem{kuijlaars2004riemann}
A.B.J. Kuijlaars, K.T.-R. McLaughlin, W.~Van~Assche, and M.~Vanlessen.
\newblock The riemann--hilbert approach to strong asymptotics for orthogonal
  polynomials on $[- 1, 1]$.
\newblock {\em Advances in Mathematics}, 188(2):337--398, 2004.

\bibitem{LevinLubinskyOPs}
A.~L. Levin and D.~S. Lubinsky.
\newblock {\em Orthogonal Polynomials for Exponential Weights}, volume~4.
\newblock Springer, 2001.

\bibitem{SOEquilibriumMeasure}
S.~Olver.
\newblock Computation of equilibrium measures.
\newblock {\em J. Approx. Theory}, 163:1185--1207, 2011.

\bibitem{SOPainleveII}
S.~Olver.
\newblock Numerical solution of {Riemann--Hilbert} problems: {Painlev\'e II}.
\newblock {\em Found. Comput. Math.}, 11:153--179, 2011.

\bibitem{SORHFramework}
S.~Olver.
\newblock A general framework for solving {Riemann--Hilbert} problems
  numerically.
\newblock {\em Numer. Math.}, 122:305--340, 2012.

\bibitem{SOTrogdonNNSD}
S.~Olver and T.~Trogdon.
\newblock Nonlinear steepest descent and the numerical solution of
  {Riemann--Hilbert} problems.
\newblock {\em Comm. Pure Appl. Maths}, 2012.
\newblock To appear.

\bibitem{SOTrogdonRMT}
S.~Olver and T.~Trogdon.
\newblock Numerical solution of {R}iemann--{H}ilbert problems: random matrix
  theory and orthogonal polynomials.
\newblock {\em Const. Approx.}, 2013.
\newblock To appear.

\bibitem{TrogdonThesis}
T.~Trogdon.
\newblock {\em Riemann--Hilbert Problems, Their Numerical Solution and the
  Computation of Nonlinear Special Functions}.
\newblock PhD thesis, Univeristy of Washington, 2013.

\end{thebibliography}

\end{document}